\documentclass[12pt]{article}
\usepackage{amssymb,latexsym,amscd,amsmath,amsfonts,amsthm,enumerate}
\usepackage{graphicx}

\newtheorem{theorem}{Theorem}[section]

\newtheorem{corollary}[theorem]{Corollary}

\newtheoremstyle{definition}
  {6pt}
  {6pt}
  {}
  {}
  {\bfseries}
  {.}
  {.5em}
  {}%
\theoremstyle{definition}
\newtheorem{definition}[theorem]{Definition}

\newtheoremstyle{example}
  {6pt}
  {6pt}
  {}
  {}
  {\bfseries}
  {.}
  {.5em}
  {}%
\theoremstyle{example}
\newtheorem{example}[theorem]{Example}

\newtheoremstyle{remark}
  {6pt}
  {6pt}
  {}
  {}
  {\bfseries}
  {.}
  {.5em}
  {}%

\theoremstyle{remark}

\makeatletter
\renewcommand\@makefntext[1]{%
\setlength\parindent{1em}%
\noindent

\makebox[1.8em][r]{}{#1}}
\makeatother

\begin{document}
\title{\bf  \large ON TANGENTS TO  CURVES }
\author{
 \centerline{\small Duong Quoc Viet}\\}
   \date{}
\maketitle \centerline{\parbox[c]{10.45 cm}{ \small {\bf
ABSTRACT:}   In this paper, we give a  simple definition of
tangents to a curve in elementary geometry.  From which, we
characterize the existence of the tangent to a curve at a point.}}

\section{Introduction}
It has long been known that the notion of tangents to a curve is
one of most important notions of analytic geometry and classical
analytic.
  The
first definition of tangents  was "a right line which touches a
curve, but which when produced, does not cut it" \cite{W}. This
old definition prevents inflection points from having any tangent.
It has been dismissed and the modern definitions are equivalent to
those of Leibniz. Pierre de Fermat developed a general technique
for determining  tangents of a curve by using his method of
adequality in the 1630s. Leibniz defined a tangent line as a line
through a pair of infinitely close points on the curve (see e.g.
\cite{T}).

The notion of  tangents to an arbitrary curve can be traced back
to the work of Archimedes in the third century B.C, when he solved
the problem of finding tangents to spirals. The geometric idea of
 tangent lines as the limit of secant lines serves as the
motivation for analytical methods that are used to find tangent
lines explicitly (see e.g. \cite{T}). The question of finding
tangent lines to a graph, or the tangent line problem, was one of
the central questions leading to the development of calculus in
the 17th century. In the second book of his Geometry \cite{D},
Ren\'e Descartes said of the problem of constructing tangents to a
curve. "And I dare say that this is not only the most useful and
most general problem in geometry that I know, but even that I have
ever desired to know"\cite{L}. \footnotetext{\begin{itemize}
\item[ ]{\it Mathematics  Subject  Classification}(2010): 53A04;
00A05. \item[ ]{\it Key words and phrases}: tangent, meaning of
derivative.
\end{itemize}}
\begin{figure}[ht]
\centering
\includegraphics{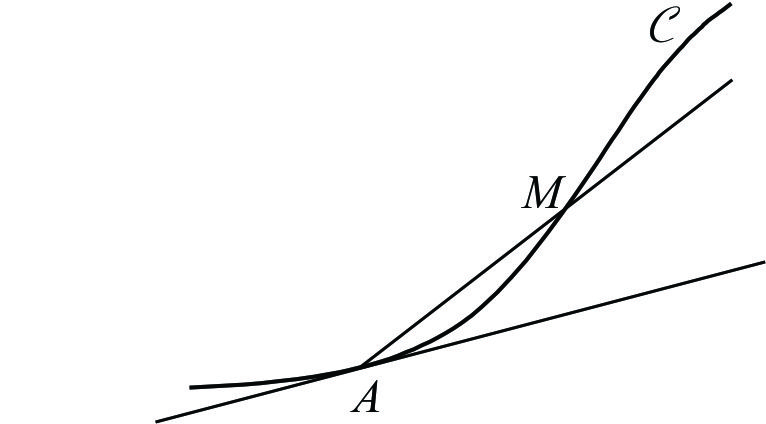}
\end{figure}

Up to now, in analytic geometry and classical analytic, one use
the following definition of tangents (see e.g. \cite{F, GB, K,
GB1}).

 \begin{definition} Consider the sequence of straight lines (secant lines) passing
through two points, $A$ and $M$, those that lie on $\cal C$. The
tangent to $\cal C$ at $A$ is defined as the limiting position of
the secant line $AM$ as $M$ tends to $A$ along the curve $\cal C$.
\end{definition}

 However,  in elementary geometry,
this definition is still very hard to explain, because what is
"the limiting position of secant lines"? And how to overcome this
problem?

In this paper, we first would like to give  a  simple definition
of tangents to a curve in elementary geometry.
 Next, we characterize the existence of  tangents to
 a curve in the space and the
relationship between the  existence of  tangents  and the
derivative.

 \section{The definition of tangents}

  In elementary geometry, it is very difficult to make the notion about "the limiting
position of secant lines" in Definition 1.1 accurate.
   Because one need to put lines
in a metric space or a topological space to consider the limit of
them.

Overcoming this problem is based on the idea of the definition of
 tangents in \cite{F} by Flett and remark that: each line
passing through a fixed point is determined by its direction
vectors. This view helps us move considering "the limiting
position of secant lines" to considering "the limiting position of
direction vectors". But the direction vectors
 of a line must be not zero.

The remark above leads us to considering  "the limiting position of direction vectors
of a  stable length".

That  is why we give the following definition.

 \begin{definition}  Let $L$  be a curve in the space and a point  $A \in L.$
 Consider an arc of $L$  containing  $A$ which is divided into two parts $L_1$ and $L_2$ by $A$
 such that $L_1$ and $L_2$ intersect at only  $A,$ and in each part there always exist points
 that are not $A$.  For each point $M \in L$ and $M \not= A$, the secant
line  $AM$  has a unit vector
$\frac{\overrightarrow{AM}}{\mid\overrightarrow{AM}\mid}.$
  Then if the limits
$$\lim_{\substack{
M \rightarrow A\\
M \in L_1
}}
\frac{\overrightarrow{AM}}{\mid\overrightarrow{AM}\mid} \; \text{and} \;
   \lim_{\substack{
M \rightarrow A\\
M \in L_2
}}
\frac{\overrightarrow{AM}}{\mid\overrightarrow{AM}\mid}
$$  exist  and they are collinear vectors,
  we call that  $L$  has the tangent at  $A,$ and both these
limits are direction vectors of this tangent.

\begin{center}
\includegraphics{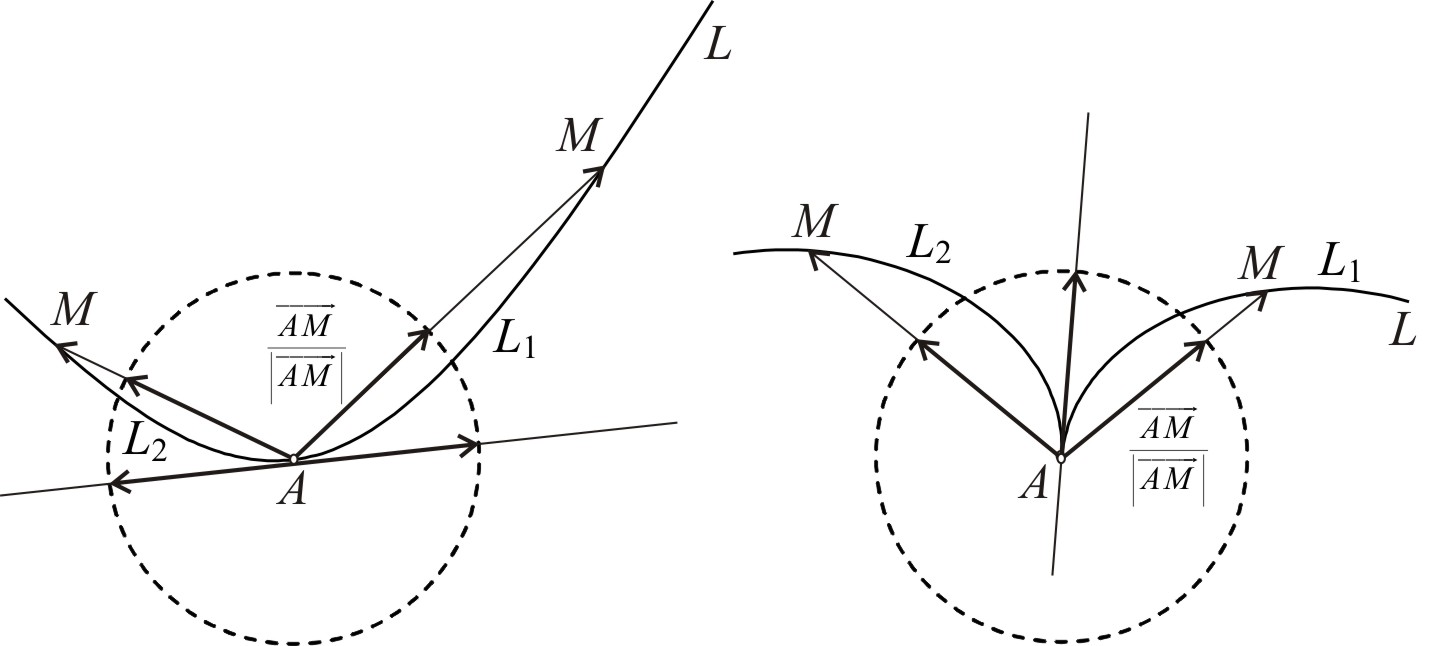}
\end{center}
\end{definition}

 \section{The condition for the existence of tangents}

 Suppose that  a curve $L$ has the  parametric equation
      $$\left \{%
\begin{array}{lll}
x &=& x(t)\\
y&=& y(t)\\
z &=& z(t)
\end{array}%
\right.$$ and  $A(x(t_0), y(t_0), x(t_0)) \in L.$  Let  $M(x(t), y(t), z(t)) \in L$, where $t \ne t_0.$  Set
$$\Delta t = t-t_0;  \Delta x = x(t)-x(t_0); \Delta y = y(t)-y(t_0); \Delta z = z(t)-z(t_0).$$
Then we have
 $$\frac{\overrightarrow{AM}}{\mid\overrightarrow{AM}\mid} =
  \frac{\Delta t}{\mid\Delta t\mid}\frac{1}{\sqrt{\big(\frac{\Delta x}{\Delta t}\big)^2 + \big(\frac{\Delta y}{\Delta t}\big)^2 + \big(\frac{\Delta z}{\Delta t}\big)^2}} \bigg( \frac{\Delta x}{\Delta t}, \frac{\Delta y}{\Delta t}, \frac{\Delta z}{\Delta t}\bigg).$$
  From this equation it follows that  $L$ has the tangent at $A$  if and only if  the
  following limits
    $$\lim_{\Delta t
 \rightarrow 0^+}\frac{1}{\sqrt{\big(\frac{\Delta x}
 {\Delta t}\big)^2 + \big(\frac{\Delta y}{\Delta t}\big)^2 +
 \big(\frac{\Delta z}{\Delta t}\big)^2}}
  \bigg( \frac{\Delta x}{\Delta t},
  \frac{\Delta y}{\Delta t}, \frac{\Delta z}{\Delta t}\bigg)$$ and $$\lim_{\Delta t
 \rightarrow 0^-}\frac{1}{\sqrt{\big(\frac{\Delta x}
 {\Delta t}\big)^2 + \big(\frac{\Delta y}{\Delta t}\big)^2 +
 \big(\frac{\Delta z}{\Delta t}\big)^2}}
  \bigg( \frac{\Delta x}{\Delta t},
  \frac{\Delta y}{\Delta t}, \frac{\Delta z}{\Delta t}\bigg) $$
exist and are collinear vectors.

   Hence we obtain the following theorem.
    \begin{theorem} Let   $L$ be a curve having the parametric equation
      $\left \{%
\begin{array}{lll}
x &=& x(t)\\
y&=& y(t)\\
z &=& z(t)
\end{array}%
\right.$ and  $A(x(t_0), y(t_0), x(t_0)) \in L.$ Set
$$\Delta t = t-t_0;  \Delta x = x(t)-x(t_0); \Delta y = y(t)-y(t_0); \Delta z = z(t)-z(t_0).$$
Then $L$ has the tangent at $A$  if and only if  the
  following limits
    $$\lim_{\Delta t
 \rightarrow 0^+}\frac{1}{\sqrt{\big(\frac{\Delta x}
 {\Delta t}\big)^2 + \big(\frac{\Delta y}{\Delta t}\big)^2 +
 \big(\frac{\Delta z}{\Delta t}\big)^2}}
  \bigg( \frac{\Delta x}{\Delta t},
  \frac{\Delta y}{\Delta t}, \frac{\Delta z}{\Delta t}\bigg)$$ and $$\lim_{\Delta t
 \rightarrow 0^-}\frac{1}{\sqrt{\big(\frac{\Delta x}
 {\Delta t}\big)^2 + \big(\frac{\Delta y}{\Delta t}\big)^2 +
 \big(\frac{\Delta z}{\Delta t}\big)^2}}
  \bigg( \frac{\Delta x}{\Delta t},
  \frac{\Delta y}{\Delta t}, \frac{\Delta z}{\Delta t}\bigg) $$
exist and these vectors are collinear.
    \end{theorem}

From this theorem, we immediately get the following result.
\begin{corollary} Let   $L$ be a curve having the parametric equation
      $\left \{%
\begin{array}{lll}
x &=& x(t)\\
y&=& y(t)\\
z &=& z(t)
\end{array}%
\right.$ and  $A(x(t_0), y(t_0), x(t_0)) \in L.$ Set
$$\Delta t = t-t_0;  \Delta x = x(t)-x(t_0); \Delta y = y(t)-y(t_0); \Delta z = z(t)-z(t_0).$$
Then $L$ has the tangent at $A$  if the
  following limit exists
    $$\lim_{\Delta t
 \rightarrow 0}\frac{1}{\sqrt{\big(\frac{\Delta x}
 {\Delta t}\big)^2 + \big(\frac{\Delta y}{\Delta t}\big)^2 +
 \big(\frac{\Delta z}{\Delta t}\big)^2}}
  \bigg( \frac{\Delta x}{\Delta t},
  \frac{\Delta y}{\Delta t}, \frac{\Delta z}{\Delta t}\bigg).$$

\end{corollary}

Particularly, if $x(t), y(t), z(t)$ are differentiable at $t_0$
and  $x'(t_0), y'(t_0), z'(t_0)$ are not all zero, then $L$ has
the tangent at $A$ and $(x'(t_0), y'(t_0), z'(t_0))$ is a
direction vector of this tangent. In  the two-dimensions  space,
if $L$ is the graph of the function $y = f(x)$ and $A \in L$ has
the abscissa $x_0$ and $f(x)$ is differentiable at  $x_0,$
  then from above result, it implies that  $L$  has the tangent at $A$  and $(1, f'(x_0))$ is
  a direction vector of this tangent, i.e., the slope of this tangent is  $f'(x_0).$

\begin{example} Let $(C)$ be a circle    with the center $I(a, b)$ and  the radius $R$.
Then it has the  parametric equation
$$\left \{%
\begin{array}{l}
x = a + R\mathrm{cos} t\\
y= b+ R\mathrm{sin} t.
\end{array}%
\right.$$  Now assume that  $$A(a + R\mathrm{cos} t_0, b+ R\mathrm{sin} t_0 ) \in (C).$$ Then the direction vector of the tangent at  $A$ is
$\overrightarrow{u} = (-\mathrm{sin} t_0,   \mathrm{cos} t_0).$  It is easily to check that
 $\overrightarrow{u}$ is perpendicular to $\overrightarrow{IA} = (R \mathrm{cos} t_0, R  \mathrm{sin} t_0).$  Therefore the tangent  to $(C)$ at  $A$ is the line passing through  $A$ and  perpendicular to the radius $IA.$  We obtain the classical  method of determining  tangents to  circles.
\end{example}

\begin{example} The curve $y = \sqrt[3]{x}$ is not differentiable at  $0,$ however since
$$\lim_{\Delta x
 \rightarrow 0}
  \frac{1}
 {\sqrt{1+\big(\frac{\sqrt[3]{\Delta x}}{\Delta x}\big)^2}}
    \bigg( 1,  \frac{\sqrt[3]{\Delta x}}{\Delta x}\bigg) = (0, 1),$$ it follows  that the
    vertical axis is the tangent to this curve at the point $O(0,0).$
\end{example}
This example shows that the problem of the existence of tangents
and the problem of  differentiability are not equivalent.

 Return to  consider closely the curve  $L$ which is the graph of the function
$y=f(x)$. Let a point $A \in L$ have the abscissa  $x_0.$  Then
$L$ has the tangent at $A$ if and only if  the following limits
$$\lim_{\Delta x
 \rightarrow 0^+}
  \frac{1}
 {\sqrt{1+\big(\frac{\Delta y}{\Delta x}\big)^2}}
    \bigg( 1,  \frac{\Delta y}{\Delta x}\bigg) \; \mathrm{and }\; \lim_{\Delta x
 \rightarrow 0^-}
  \frac{1}
 {\sqrt{1+\big(\frac{\Delta y}{\Delta x}\big)^2}}
    \bigg( 1,  \frac{\Delta y}{\Delta x}\bigg)$$
     exist and these vectors are collinear vectors.  But it is clear that this condition is equivalent to  that
$\lim_{\Delta x \rightarrow 0} \frac{\Delta y}{\Delta x}$ is
finite or infinite. Concretely, if $$\lim_{\Delta x \rightarrow
0^+} \frac{\Delta y}{\Delta x}= \pm \infty \; \mathrm{and}\;
\lim_{\Delta x \rightarrow 0^-} \frac{\Delta y}{\Delta x}= \pm
\infty,$$
$$\lim_{\Delta x
 \rightarrow 0^+}
  \frac{1}
 {\sqrt{1+\big(\frac{\Delta y}{\Delta x}\big)^2}}
    \bigg( 1,  \frac{\Delta y}{\Delta x}\bigg)= (0,\pm 1) \; \mathrm{and}\; \lim_{\Delta x
 \rightarrow 0^-}
  \frac{1}
 {\sqrt{1+\big(\frac{\Delta y}{\Delta x}\big)^2}}
    \bigg( 1,  \frac{\Delta y}{\Delta x}\bigg)= (0,\pm 1)$$ and then the tangent to
    $L$ at  $A$ is parallel to the vertical axis, i.e., the slope of the tangent is infinite.
If $\lim_{\Delta x \rightarrow 0} \frac{\Delta y}{\Delta x}$ is finite, i.e., $f(x)$ is  differentiable at $x_0$ and
    $\lim_{\Delta x \rightarrow 0} \frac{\Delta y}{\Delta x}= f'(x_0), $
    then $$\lim_{\Delta x
 \rightarrow 0}
  \frac{1}
 {\sqrt{1+\big(\frac{\Delta y}{\Delta x}\big)^2}}
    \bigg( 1,  \frac{\Delta y}{\Delta x}\bigg)= \frac{1}
 {\sqrt{1+(f'(x_0))^2}}
    \bigg( 1,  f'(x_0)\bigg),$$ and the slope of the tangent at $A$ is $f'(x_0).$

To express the relationship between the  existence of tangents
and the derivative, we consider
 $\lim_{\Delta x \rightarrow 0} \frac{\Delta y}{\Delta x}= \infty$ as the derivative of
 the function at  $x_0$ and write \linebreak
 $f'(x_0) = \infty.$ And we call   {\it extended derivative} of the function $y= f(x)$
 at $x_0$ the  finite or  infinite limit  of
$\lim_{\Delta x \rightarrow 0} \frac{\Delta y}{\Delta x}.$ Then we
obtain the following corollary.

\begin{corollary} The graph of the function
$y= f(x)$ has the tangent at $A $ with the  abscissa $x_0$ if and only if  $f(x)$ has the extended derivative at  $x_0,$ and in this case the slope of the tangent is $f'(x_0).$
\end{corollary}

\noindent{\bf Acknowledgement:} Thanks are due also to my
colleagues, several  members  of Department of Mathematics of
Hanoi National University of Education, and special thanks are due
to Dr. Truong Thi Hong Thanh for  excellent advice and support.

\noindent
Department of Mathematics\\
\noindent
Hanoi National University of Education\\
\noindent
136 Xuan Thuy street, Hanoi, Vietnam\\
\noindent
 vietdq@hnue.edu.vn


\begin{thebibliography}{99}
\bibitem {D} R.  Descartes, {\it The geometry of Ren\'e
Descartes}, Courier Dover. pp. 95. ISBN 0-486-60068-8(1954).

\bibitem {F}T. M. Flett, {\it The Definition of a Tangent to a Curve}, Edin.
Math. Notes\;\;{\bf l41}(1957), 1-9. DOI:
10.1017/S0950184300003153, Published online: 31 October 2008.

\bibitem {GB}S. R. Ghorpade, V. Balmohan, {\it A Course in Multivariable Calculus and Analysis},
   Springer Science + Business Media,  LLC, 2006.
\bibitem {K}  K.  E. Hirst, {\it Calculus of One Variable}, Springer-Verlag London Limited, 2006.

\bibitem {GB1} W. Kühnel, {\it Differential Geometry: Curves-Surfaces-Manifolds},
Amer. Math. Soc, 2006.
\bibitem {L}R. E. Langer, {\it Ren\'e Descartes}, American
Mathematical Monthly (Mathematical Association of America) 44 (8):
495-512. doi:10.2307/2301226. JSTOR 2301226 (October 1937).
\bibitem {W}  N.  Webster, {\it American Dictionary of the English
Language}, vol.2, p.733  (New York: S. Converse, 1828).
\bibitem {T} en.wikipedia.org/wiki/Tangent.

\end{thebibliography}
\end{document}